\documentclass[12pt]{article}
\begin{document}
 \def\K{I\!\!K}
\def\N{\mathbb N}
\def\Z{\mathbb Z}
\def\A{\mathbb A}
\def\B{\mathbb B}
\def\C{\mathbb C}
\def\D{\mathbb D}
\def\E{\mathbb E}
\def\R{I\!\!R}
\def\H{I\!\!H}
\def\C{I\!\!\!\!C}
\def\P{I\!\!P}
\def\Q{I\!\!Q}
\def\B{I\!\!B}
\def\F{I\!\!F}
\def\M{I\!\!M}
\def\D{I\!\!D}
\def\P{I\!\!P}
\def\r{\rho}
\def\l{\lambda}
\def\s{\sigma}
\def\b{\beta}
\def\a{\alpha}
\def\d{\delta}
\def\g{\gamma}
\def\z{\zeta}
\def\fp{\hfill \Box \\[2pt]}
\def\e {\epsilon}
\def\proj{\mathop {\rm proj}}
\def\ind{\mathop {\rm ind}}
\def\lim{\mathop {\rm lim}}
\def\id{\mathop {\rm id}}
\def\ker{\mathop {\rm ker}}
\def\max{\mathop {\rm max}}
\def\min{\mathop {\rm min}}
\def\sup{\mathop {\rm sup}}
\def\supp{\mathop {\rm supp}}
\textwidth 6.01in\textheight 9in \evensidemargin 0pt
\oddsidemargin 0pt\topskip -0in\topmargin 0pt
\begin{center}
\Large{\bf  LAPLACE EQUATIONS FOR THE DIRAC, EULER OPERATORS AND THE HARMONIC OSCILLATOR\\
Ahmedou Yahya ould Mohameden  and  Mohamed Vall ould Moustapha}
\end{center}
{\bf ABSTRACT:} In this article, we give the explicit solutions to the Laplace equations associated to the Dirac operator, Euler operator and  Harmonic oscillator in $\R$.
\begin{center}
{\bf 1 -- Introduction:}
\end{center}

In the paper $[4]$ we have solved the wave Cauchy problems associated to the Dirac, Euler operator and the harmonic oscillator.\\
The aim of this paper is to give the explicit solutions to the following Cauchy problems of Laplace type
$$ \left \{\begin{array}{cc}\left(D+\frac{\partial^2}{\partial y^2}\right)u(y,X)=0&;(y,X)\in \R^\ast_+\times \R\\ u(0,X)=u_0(X) ~~; u_0\in C^\infty_0(\R)\end{array}
\right.\eqno(LD) $$
$$ \left \{\begin{array}{cc}\left(E_{a}+\frac{\partial^2}{\partial y^2}\right)v(y,\xi)=0&;(y,\xi)\in \R^\ast_+\times \R\\ v(0,\xi)=v_{0}(\xi) ~~; v_0\in C_0^\infty(\R)\end{array}
\right. \eqno(LE_{a})$$
$$ \left \{\begin{array}{cc}\left(H^{a}+\frac{\partial^2}{\partial y^2}\right)w(y,x)=0;(y,x)\in \R^\ast_+\times \R\\ w(0,x)=w_{0}(x) ~~; w_0\in
C_0^\infty(\R)\end{array}
\right. \eqno(LH_{a})$$
where:
$$D=\frac{\partial}{\partial X},\eqno(1.1)$$
$$E^{a}=-2a\xi\frac{\partial}{\partial\xi},\ \ \ \ \ \ \ \  a> 0 \eqno(1.2)$$
$$H^a=\frac{\partial ^2}{\partial x^2} - a^2 x^2,\ \ \ \ \ \ \ a> 0 \eqno(1.3)$$
are respectively
the Dirac operator, the Euler operator and the harmonic oscillator on $\R$.
 These operators
  play a fundamental role in many mathematical and physical problems. In physics the harmonic oscillator appears when e.g. modeling atoms and their
quantum states $\left([2],[3]\right)$.
Note that the Laplace equation associated to the Euler operator $(LD)$ is a typical example of a hyperbolic equation with multiple characteristics considered by Leray in 1960 $\left( \mbox{\rm see} [1], P.355\right)$.
\begin{center}
{\bf 2-- Laplace equation for the Dirac and the Euler operators }
\end{center}
Here our objective is to solve the Cauchy problems $(LD)$ and $(LE_a)$, for this we need the following recall:\\
Let $\Omega\subset \R$
be an open set and let $\nu \geq 1$ be a fixed real number.\\
The class $G^\nu(\Omega)$ of Gevrey functions of order $\nu$ in $\Omega $ is the set of functions $ f \in C^\infty(\Omega)$
satisfying the property that for every compact subset $K$ of $\Omega$, there exists a positive
constant $C = C_K$ such that for all $ l \in N $ and all
$x\in K$
$$|\partial^lf(x)| \leq C^{
l+1}(l!)^\nu.$$
 It is easy to recognize
that $G^1= A(\Omega)$, the space of all analytic functions in $\Omega$.\\
Assume $\nu> 1$, we shall denote by $G^\nu_0(\Omega)$
the vector space of all $f \in G^\nu(\Omega)$
with compact
support in $\Omega$ $\left(\mbox{\rm see}[5]\right)$.\\
{\bf Theorem 2.1:} The Cauchy problem $(LD)$ for the Laplace equation associated to the Dirac operator has the unique solution given by:\\
$$U(y,X) = \int_{X}^{+\infty}P_{D}\left(y,X,X'\right)U_0(X')dX'\eqno(2.1)$$
where $$P_{D}(y,X,X')=\frac{y}{2\sqrt{\pi}}\left(X'-X\right)^{-3/2}\exp\left\{\frac{-y^{2}}{4(X'-X)}\right\}\eqno(2.2)$$
{\bf Proof}.  Note that the Cauchy Problem for the Laplace equation associated to the Dirac operator
is well-posed in $G^\nu(\R)$ for any $1<\nu < 2$  see $[7]$.\\
The fact that the function  $P_{D}(y,X,X')$ satisfies the Laplace equation can be checked by direct computation.
To see the limit condition we use the change of variables  $s=\frac{y^2}{4(X'-X)}$.\\
{\bf Corollary 2.2:} The Cauchy problem $(LE_{a})$ for the Laplace equation associated to the Euler operator has the unique solution given by:\\
$$V(y,\xi)~=~\int_{|\xi'|<|\xi|}P_{E^{a}}(y,\xi,\xi')V_{0}(\xi')d\xi'\eqno(2.3)$$
where $L_{E^{a}}(y,\xi,\xi')$ is given by
 $$P_{E^{a}}(y,\xi,\xi')=~\sqrt{\frac{a}{2\pi}}\frac{y}{\xi'}log^{-3/2}(|\xi/\xi'|)\exp\left\{\frac{-ay^{2}}{2log(|\xi/\xi'|)}\right\} \ \eqno(2.4)$$
{\bf Proof:} This can be shown by using the change of variable $X=\frac{log |\xi|}{-2 a}$ and the theorem 2.1.\\
\begin{center}{\bf 3--Laplace equation for the Harmonic oscillator}
\end{center}
{\bf Theorem 3.1} The Cauchy problem $(LH^{a})$ for the Laplace equation associated to the Harmonic oscillator has the unique solution given by:
$$w(y,x) = \int_{-\infty}^{+\infty}P_{H^{a}}\left(y,x,x'\right)w_{0}(x')dx'\eqno(3.1)$$
where $P_{H^{a}}(y,x,x')$ is given by\\
 $$P_{H^{a}}\left(y,x,x'\right)=\frac{\sqrt{a}y}{2\pi}\int^{+\infty}_{0}u^{-3/2}.\left(sinh(2a u)\right)^{-1/2}\times $$
 $$exp\left\{\frac{-y^{2}}{4 u}-\frac{a}{2}\left(x^{2}+x'^{2}\right)coth(2 a u)+\frac{axx'}{sinh(2 a u)}\right\}d u\eqno(3.2)$$
{\bf Proof:}  We recall the integral
kernel of the heat operator for the harmonic oscillator (Mehler-Fuchs- formula, $[2]$)
 $$K_{a}(t,x,x')=\sqrt{\frac{a}{2\pi}}\frac{1}{\sqrt{\sinh(2at)}}\times$$
 $$~exp\left[-\frac{a}{2}(x^{2}+x'^{2})\coth(2at)+\frac{axx'}
 {\sinh(2at)}\right]\eqno(3.3)$$
Taking the derivative of the formula $\left(6],p.50\right)$
$$\frac{e^{-t \lambda}}{t}=\frac{1}{\sqrt{\pi}}\int_0^\infty e^{-u t^2}u^{-1/2}e^{-\lambda^2/4u}du\eqno(3.4)$$
 with respect to $\lambda$ we can write
$$\frac{e^{-t \lambda}}{t}=\frac{1}{2\sqrt{\pi}}\int_0^\infty e^{-u t^2}u^{-3/2}\lambda e^{-\lambda^2/4u}.du\eqno(3.5)$$
By setting $t=\sqrt{-H_a}$ and $\lambda=y$ in $(3.5)$ we get
$$e^{-y\sqrt{-H_a}}=\frac{y}{\sqrt{2\pi}}\int_0^\infty e^{-s^2/4u}u^{-3/2}e^{u H_a}du\eqno(3.6)$$
where $ e^{u H_a}$ is the heat operator of the harmonic oscillator.
Let $P_{H_a}\left(y,x,x'\right)$ be the integral kernel of $e^{-y\sqrt{-H_a}}$, making use of the formula $(3.3)$ in $(3.6)$ we see that the formula $(3.2)$ holds  and the proof of the theorem 3.1 is finished.\\
{\bf Proposition 3.1} Let $P(y, x,x') = \frac{1}{\pi}\frac{y}{y^{2}+(x-x')^{2}}$ be the Poisson kernel on the half- plane $\left([6],p.60\right)$. Then we have
$$lim_{a\searrow 0}P_{H_{a}}(y,x,x')~=~P(y, x,x').\eqno(3.7)$$
{\bf Proof:} The formula $(3.7)$ uses essentially the formula $(3.2)$ and $lim_{s\searrow 0}~\frac{s}{sinh(s)}~=~1$.
\begin{center}{\bf
4-Directions for further studies:}
\end{center}
We suggest here a certain number of open related problems connected to this paper.
We are interested in the Laplace equation for the harmonic oscillator with an inverse square potential
$$ \left \{\begin{array}{cc}\frac{\partial^2}{\partial y^2}
u(y,x)+\left(\frac{\partial^2}{\partial x^2} - a^2 x^2-\frac{b^2}{x^2}\right)u(y,x)&(y,x)\in \R_+^\ast\times \R\\ ~u(0,x)=f(x) & f\in
C_0^\infty(\R)\end{array}
\right.$$\\
Finally, we suggest problems in direction of the non linear Laplace equations for the harmonic oscillator
and to look for global solution and a possible  blow up in finite times.

\begin{flushleft}

Universit\'e de Nouakchott Al-asriya\\
Facult\'e des Sciences et Techniques\\
D\'epartement de Math\'ematiques et Informatique\\
Unit\'e de Recherche: {\bf \it Analyse, EDP et Mod\'elisation: (AEDPM)}\\
B.P: 5026, Nouakchott-Mauritanie\\
E-mail address: ahmeddou2011@yahoo.fr\\

Al Jouf University\\
College of Sciences and Arts, Al-Qurayyat
Saoudi Arabia.\\
Universit\'e de Nouakchott Al-asriya\\
Facult\'e des Sciences et Techniques\\
D\'epartement de Math\'ematiques et Informatique\\
Unit\'e de Recherche: {\bf \it Analyse, EDP et Mod\'elisation: (AEDPM)}\\
B.P: 5026, Nouakchott-Mauritanie\\
E-mail adresse:mohamedvall.ouldmoustapha230@gmail.com\\
\end{flushleft}

\begin{thebibliography}{999}
\bibitem{1} Armand Borel, Gennadi M. Henkin, and Peter D. Lax; JEANT LERAY
(1906---1998) NOTICES OF THE AMS VOLUME 47, NUMBER 3 (2000).
\bibitem{1} Berline, N.,Getzler, E., Vergne, M. Heat kernels and dirac operator Springer Verlag 2004.
%\bibitem{2} Colombini, F. and Nishitani, T.
%On finitely finitely degenerate hyperbolic operators of second order, Osaka J. Math.
%41 (2004), 933�947
\bibitem{3} G. B. Folland, Quantum Field Theory, A Tourist Guide for Mathematicians ; Mathematical surveys and Monographs vol 149 (2008)
%\bibitem{4} Greiner P. C. , Daniel Holocman, and Yakar Kannai, Wave kernels related to the second order operator, Duke Math.J.vol. 114,  329-387 (2002).
%\bibitem{5} Guelfand, I.M., and Chilov, C. E. Les distributions Tome 3 Th\'eorie des \'equations %diff\'erentielles, Dunod 1968.
    %\bibitem{6}  Y.Kannai, The method of ascent and $\cos\sqrt{A^2+B^2}$
%Bull. Sci. Math. 124 (2000),
%573 -597. MR 2001h:35025 333.
%\bibitem{7}
%Magnus, W. Oberhittenger, F. and Soni R. P.  Formulas and Theorems for the special functions %of Mathematical
   % physics,Springer-Verlog New-York 1966.
%\bibitem{8} Dott. Marcello D�Abbicco, On the Cauchy Problem
%for Linear Hyperbolic
%Equations and Systems,Tesi di Dottorato di Ricerca in Matematica
%Universit`A Degli Studi dI Bari
%XX Ciclo � A.A. 2007�2008
%Settore Scientifico: MAT/05
\bibitem{2}Ould Mohameden A. Y. and Ould Moustapha M. V.;
 Wave kernels of the Dirac, Euler operators and the harmonic oscillator (submitted to JMP)\\
 \bibitem{9}   L. Rodino. Linear partial differential operators in Gevrey spaces. World Scientific
Publishing Co. Inc., River Edge, NJ, 1993.
\bibitem{10}  Strichartz Robert S. A guide to distribution theory and Fourier transform, Studies in advanced mathematics
CRC press, Boca racon Ann Arbor london tokyo 1993.
%\bibitem{[11} M.E. Taylor, Partial Differential Equations I, Springer, New York, Berlin, %Heidelberg, 1996.
\bibitem{10} Tamotu Kinoshita and Giovanni Taglialatela;
Time regularity of the solutions to second order
hyperbolic equations,



Ark. Mat., 49 (2011), 109–127
\end{thebibliography}
\end{document}